\documentclass[12pt]{article}
\usepackage{amsthm}
\usepackage{amsmath,amsfonts,amssymb,txfonts,mathrsfs}

\usepackage[all]{xy}
\newtheorem{thm}{Theorem}[section]
\newtheorem{lem}[thm]{Lemma}
\newtheorem{proposition}[thm]{Proposition}
\newtheorem{defin}{Definition}
\newtheorem{corol}[thm]{Corollary}

\theoremstyle{remark}

\author{Saurav Bhaumik and Vikram Mehta}
\title{Harder-Narasimhan Filtrations which are not split by the Frobenius maps}
\date{}
\newcommand{\spec}{{\rm Spec\;}}

\begin{document}
\maketitle

\section{Introduction} Let $X$ be a smooth projective variety over a perfect field $k$ of characteristic $p>0$, and $V$ be a vector bundle over $X$. Recall that in characteristic $p$, a vector bundle $W$ is called \emph{strongly semistable} if $(F^n)^*(W)$ is semistable for all $n$, where $F:X\to X$ is the absolute Frobenius. If $X$ is a curve and $V$ is not strongly semistable, then some Frobenius pullback $(F^t)^*V$ is a direct sum of strongly semistable bundles (Proposition 2.1, \cite{biswasparam} or Corollary 5.2, \cite{langersemistableg}). A natural question to ask is whether this still holds for ${\rm dim}\, X>1$. Biswas, et al. in \cite{biswas} showed that there is always a counterexample to this over any algebraically closed field of positive characteristic which is uncountable. However, we will produce a smooth projective variety over $\mathbb Z$ and a rank $2$ vector bundle on it, which, restricted to each prime $p$ in a nonempty open subset of $\spec\mathbb Z$, constitutes a counterexample over $p$. Indeed, given any split semisimple simply connected algebraic group $G$ of semisimple rank $>1$ over $\mathbb Z$, we will show that there exists a smooth projective homogeneous space $X_Z$ over $\mathbb Z$ and a vector bundle $V$ on $X_Z$ of rank $2$ such that for each prime $p$ in a nonempty open subset of $\spec\mathbb Z$, the restriction $V\otimes\mathbb F_p$ as a vector bundle over $X_Z\otimes\mathbb F_p$ is a counterexample. We only use the Borel-Weil-Bott theorem in characteristic 0 (see Borel \cite{borelarcata}) and Frobenius Splitting of $G/B$ in characteristic $p$ (see Mehta-Ramanathan \cite{mehtaramanathanannals}).

Let $X_Z=G/B$, where $G$ is a split semisimple simply connected algebraic group over $\mathbb Z$ with semisimple rank $>1$, and $B$ a Borel subgroup containing a split maximal torus $T$. This $X_Z$ is defined over $\mathbb Z$, let $X_Q$ denote its pullback to $\spec \mathbb Q$. We will produce a line bundle $L_Z$ over $X_Z$ such that if $L_Q$ denotes the corresponding line bundle on $X_Z$, then $H^1(X_Q,L_Q)\ne 0$. Since the base change $\spec\mathbb Q\to\spec\mathbb Z$ is flat, we have \[H^1(X_Z,L_Z)\otimes \mathbb Q=H^1(X_Q,L_Q).\]Therefore $H^1(X_Z,L_Z)\ne 0$. By semicontinuity, $H^1(X_{p},L_p)\ne 0$ where $X_p$ is the fibre over $\mathbb Z/p\mathbb Z$ for any prime $p$, since the subset of $\spec \mathbb Z$ where the first cohomology vanishes is a closed set and contains the generic point given by $\mathbb Q$. Again, since $X_Z$ is projective and thus $H^1(X_Z,L_Z)$ is finite $\mathbb Z$-module, it must have a torsion free part. Take a basis of the free part, and choose one element from the basis. This element, which we will call $\theta\in H^1(X_Z,L_Z)$, remains nonzero under the natural map $H^1(X_Z,L_Z)\otimes\mathbb F_p\to H^1(X_p,L_p)$ for all primes $p$ lying in a nonempty open subset of $\spec\mathbb Z$ (because, for example, by semi-continuity, the dimension of $H^1(X_p,L_p)$ is constant on an open subset of $\spec\mathbb Z$,  and therefore the natural map is an isomorphism on that open set). Under the identification $H^1(X_Z,L_Z)={\rm Ext}^1(\mathscr O_{X_Z},L_Z)$, let the element $\theta$ denote the extension \[0\to L_Z\to V\to \mathscr O\to 0.\]Since for each $p$, the restriction $X_p$ of our scheme $X_Z$ is Frobenius split, if $M$ is any quasicoherent sheaf, the natural map $H^i(X_p,M)\to H^i(X_p,F_{p*}F_p^*M)=H^i(X,F_p^*M)$ is actually an injection, where we denote by $F_p$ the absolute Frobenius of $X_p$. Therefore, the image of the element $\theta$ in $H^1(X_p,L_p)$ remains nonzero after successive application of $F_p$, and the exact sequence it stands for is thus always non-split even after successive application of $F_p^*$.

We will select a very ample line bundle $H$ on $X_Z$ in such a way that this $L$ will have degree positive, so that $V$ is {\em not} semistable at all; since the quotient of $V$ by $L_Z$ is $\mathscr O_{X_Z}$, we see that $0\to L_p\to V_p$ remains the Harder-Narasimhan filtration after successive pullbacks by $F_p$. Since this is not split, this will produce our desired example. 

\section{Construction of the Line Bundle $L$} We will apply Borel-Weil-Bott theorem. Let $G$ be a split semisimple simply connected algebraic group over $\mathbb Z$, $B$ a Borel subgroup containing a split maximal torus $T$ of rank $n>1$. We define $X_Z=G/B$. Therefore, one always has a dominant weight $\lambda_0$ with respect to a fixed basis of simple roots, and $\rho$ being the half sum of positive roots, $w$ being a length one element of the Weyl group, $s:=w(\lambda_0+\rho)-\rho$ has the property that $w(s+\rho)-\rho$ is dominant. Hence, by Borel-Weil-Bott theorem, we have the line bundle $L_Q$ defined over $X_Q$ such that $H^1(X_Q,L_Q)\ne 0$. Note that this $L_Q$ comes from a character of a split maximal torus of $G$, so that $L_Q$ is also defined over $\mathbb Z$ i.e. there is a line bundle $L_Z$ defined over $X_Z$ whose restriction to $X_Q$ is $L_Q$.

\section{Selection of The Very Ample Line Bundle $H$:}
  Given any character $\lambda$ of the maximal torus $T$, we have the line bundle $L_\lambda$ on $X_Z$, and this gives an isomorphism $X^*(T)\cong {\rm Pic}(X_Z)$. Then ${\rm Pic}(X_Z)$ admits a basis of line bundles coming from the simple roots $\omega_1,\ldots,\omega_n$; these simple roots (and therefore the basis of $Pic(X_Z)$) has the property that a line bundle on $X_Z$ is very ample if and only if it is a strictly positive linear combination of the basis elements. Note that Kempf's vanishing theorem (Theorem 3.1, \cite{haboush1}) is independent of the particular choice of simple roots (if $\Delta$ is a base for roots, then so is $-\Delta$). Therefore the argument in section 1 shows that if $L=\sum_im_i\omega_i$, and $H^1(X_Q,L_Q)\ne 0$, all of $m_i$ are nonzero and there is $0<r<n$ such that, by suitable reordering of the indices, $m_1,\ldots,m_r$ are all positive, and the rest are negative. 

Let $N_1,\ldots,N_r$ be positive integers. Then the line bundle $H=\sum_{i\le r}N_i\omega_i+\sum_{i>r}\omega_i$ is very ample. The degree of $L$ with respect to this very ample line bundle $H$ is $L\cdot H^{d-1}$, where $d={\rm dim} X_Z$. This can be made positive by choosing the $N_i$ to be large. Indeed, all $\omega_i$ are generated by global sections, so that any monic monomial of dimension zero in these $\omega_i$ has non-negative degree.

\vskip 1em
\noindent{\bf Acknowledgements.} The second author would like to thank the Tata Institute for its hospitality.

\vskip 2em
{\it 
\noindent Addresses:
\vskip 0.5em
\noindent School of Mathematics, Tata Institute of Fundamental Research,
Homi Bhabha Road, Mumbai 400005, India. e-mail:
saurav@math.tifr.res.in
\vskip 0.5em
\noindent 203 Walkeshwar Road, Mumbai 400006, India. e-mail:
vikram@math.tifr.res.in\\ vmehta729@gmail.com}
\end{document}